\documentclass[12pt]{article}
\usepackage[vmargin=1in,hmargin=1in]{geometry}

\usepackage{amssymb, cite}
\usepackage{graphicx}

\newcommand{\be}{\begin{equation}}
\newcommand{\ee}{\end{equation}}
\newcommand{\bqn}{\begin{eqnarray}}

\newcommand{\eqn}{\end{eqnarray}}

\newcommand{\bd}{\begin{description}}
\newcommand{\ed}{\end{description}}
\newtheorem{Theorem}{Theorem}[section]

\newtheorem{remark}[Theorem]{Remark}
\newtheorem{notation}[Theorem]{Notation}
\newtheorem{claim}[Theorem]{Claim}
\newtheorem{corollary}[Theorem]{Corollary}

\newtheorem{stat}{}[section]
\newtheorem{definition}[Theorem]{Definition}
\def\bs{\begin{stat}}
\def\es{\end{stat}}
\def\ben{\begin{enumerate}}
\def\een{\end{enumerate}}

\def\bp{\noindent{\bf Proof}  \ }
\newcommand{\ep}{\hfill $\square$}

\makeatletter
\@addtoreset{equation}{section}

\makeatother

\begin{document}

\begin{center}
{\large {\bf ON GRAPHS UNIQUELY DEFINED 
\\[2ex]
BY THEIR $K$-CIRCULAR MATROIDS}}
\\[5ex]
{\large {\bf Jos\'e F. De Jes\'us}}
\\[2ex]
{\bf University of Puerto Rico, San Juan, Puerto Rico, United States}
\\[4ex]
{\large {\bf Alexander Kelmans}}
\\[2ex]
{\bf University of Puerto Rico, San Juan, Puerto Rico, United States}
\\[0.5ex]
{\bf Rutgers University, New Brunswick, New Jersey, United States}

\end{center}

\date{}

\vskip 3ex

\begin{abstract}

In 30's Hassler Whitney considered and completely solved the problem 
$(WP)$ of describing the classes of graphs $G$ having the same cycle matroid $M(G)$ 
\cite{W2,W3}. 
A natural analog $(WP)'$ of Whitney's problem  $(WP)$
is  to describe the classes of graphs $G$ having the same  matroid $M'(G)$, where $M'(G)$ is a matroid (on the edge set of $G$) distinct from $M(G)$.
For example, the corresponding problem $(WP)' = (WP)_{\theta }$ 
for 
the so-called bicircular matroid $M_{\theta }(G)$ of graph $G$
was solved in \cite{CGW,Wag}.
 In \cite {KK1} we introduced and studied the so-called $k$-circular matroids $M_k(G)$ for every non-negative integer $k$  that is a natural generalization  of the cycle matroid $M(G):= M_0(G)$ and of the bicircular matroid 
 $M_{\theta }(G):= M_1(G)$ of graph $G$. 
In this paper (which is a continuation of our paper \cite {KK1}) we 
establish some properties of graphs guaranteeing that the graphs are uniquely defined by their $k$-circular matroids.


 \vskip 1ex


\end{abstract}

\section{Introduction}
\label{Intro}

\indent

In 30's Hassler Whitney  developed a remarkable theory on the matroid isomorphism 
and the matroid duality of graphs 
\cite{W1,W2,W3,W4}.
He considered a graph $G$ and the so called cycle matroid $M(G)$ of 
$G$ 
  (whose circuits  are the  edge subsets of the cycles in $G$)
  and stated the following natural problems on pairs $\langle G, M(G)\rangle$:

\vskip 1ex
 
$(WP)$ describe the classes  of graphs having the same cycle matroid and, in particular, graphs that can be reconstructed
from cycle matroid (up to the names  of vertices) and

\vskip 1ex

$(WP^*)$
describe the pairs of graphs whose cycle matroids are dual, i.e. describe the class of graphs closed under their cycle matroids duality.

\vskip 1.5ex

Classical Whitney's graph matroid-isomorphism theorem and 
Whitney's planarity criterion provide the answers to the above questions \cite{W2,W3,W4} (see also \cite{Ox}). 

Naturally, Whitney's problems and interesting results along this line prompted further questions and research on possible strengthenings as well as various extensions or analogs of some Whitney's results
(see, for example, \cite{CGW,HalJ,HJK,Ksisd,K3skHng,Wag }).

\vskip 1.5ex

In \cite {KK1} we introduced and studied the so-called $k$-circular matroids $M_k(G)$ 
of graph $G$. It is natural to consider an analog $(WP)_k$ 
of Whitney's problem 
$(WP):= (WP)_0$ on the classes  of graphs having the same $k$-circular matroid. In this paper (which is a continuation of our paper \cite {KK1}) we consider a particular problem  of $(WP)_k$ on graphs uniquely defined by their $k$-circular matroids.

\vskip 1.5ex

Section  \ref{Notations} provides 
some notions, notation, and some necessary preliminary facts on 
matroids and  graphs. In particular, we introduce a classification of the vertex stars of a graph $G$ in terms of the corank of matroid $M_k(G)$. This classification will play a key role in the study of the problem  
on describing the graphs uniquely defined by their $k$-circular matroids.

\vskip 1.5ex

In Section \ref{Stars-Cocircuits}  we describe the vertex stars in a graph $G$ that are cocircuits of matroid $M_k(G)$.

\vskip 1.5ex

In Section  \ref{NonSep}  we characterize the so-called non-separating cocircuits of matroid $M_k(G)$ in terms of graph $G$. This notion will be essential in proving some results on pairs $\langle G, M_k(G)\rangle$ analogous to Whitney's matroid isomorphism theorem.

\vskip 1.5ex

By Whitney's matroid isomorphism theorem a graph with at least 4 vertices is uniquely defined by its cycle matroid if and only if the graph is multi $3$-connected. 
In Section \ref{UniquelyM1} we provide an extension of above Whitney's result by describing a class of graphs uniquely defined by their $k$-circular matroid. Every  vertex star in a graph $G$ of this class is a non-separating cocircuit of $M_k(G)$.
This result as well as the main result of the next section is based on the fact that  if $S$ is a non-separating cocircuit of $M_k(G)$, then $S$ is not only a vertex star of $G$ but also a vertex star of every graph $G'$ with $M_k(G')=M_k(G)$.

\vskip 1.5ex

In Section \ref{UniquelyM2} we describe another class of graphs uniquely defined by their $k$-circular matroid.  
For each graph $G$ in this class every vertex star of $G$ except for one is a non-separating cocircuit of $M_k(G)$.

\section{Main notions and notation}
\label{Notations}

\indent

The basic notions and notation 
we use here is the same as in our paper \cite{KK1} since this paper is a continuation of  \cite{KK1} (see also \cite{BM,D,Ox,Welsh}). We will remind here some of these notions and notation. 

 \vskip 1.5ex

A graph $G$ is called {\em cacti-graph} if G has no isolated vertices, no leaves, and  no cycle components. A connected cacti-graph is called a 
{\em cactus}. Let ${\cal G}_{\bowtie}$ denote the set of cacti-graphs and 
${\cal CG}_{\bowtie}$ denote the set of connected graphs from ${\cal G}_{\bowtie}$, and so each member of ${\cal CG}_{\bowtie}$ is a cactus.

 \vskip 1.5ex
 
 Let ${\cal G}$ be the set of finite graphs and $G \in {\cal G}$. 
Let $\Delta (G) = |E(G)| - |V(G)|$. 
\vskip 0.3ex
\noindent
 Instead of $\Delta (G)$ we will write simply $\Delta G$. Let $X \subseteq E(G)$  
and $G\langle X \rangle $ be the subgraph of $G$
\vskip 0.3ex
\noindent
 induced by $X$. Then $\Delta (G\langle  X \rangle)  = |X| - |V(G\langle X \rangle )|$. 

\vskip 1.5ex

For a graph $G = (V, E, \phi )$ and  
$k \ge 1$, let 
${\cal C}_{k}(G) = \{ C \subseteq E: \Delta G \langle C \rangle = k ~and ~
G \langle C \rangle \in {\cal G}_{\bowtie}$\}.   
Then ${\cal C}_{k}(G)$
is the collection of circuits of a matroid on $E$ (see \cite{KK1}). We call 

\vskip 0.3ex
\noindent
$M_k(G)=(E, {\cal C}_{k}(G))$ the  {\em $k$-circular matroid of graph 
$G$}. 

\vskip 1.5ex
Let  ${\cal I}_k(G)$, ${\cal B}_k(G)$, and ${\cal C}^*_k(G)$ denote the families of  independent sets, bases, and cocircuits
\vskip 0.3ex
\noindent
 of $M_k(G)$, respectively. 
Let 
$\rho _k(G)$ and $\rho  ^*_k(G)$ denote the rank and the corank of matroid 
\vskip 0.3ex
\noindent
$M_k(G)$.

\vskip 1.5ex

For a graph $G$ with at least one cycle, let 
$\lfloor G \rfloor$ denote the maximum subgraph of $G$ with no leaves and no isolated vertices.  Graph  $\lfloor G \rfloor$ is called the {\em kernel of graph $G$}. If $F$ is a forest, then the kernel of $F$ is not defined.

  \vskip 1.5ex
  
For a graph $G$ having a component with at least two cycles, let $ [ G] $ denote the maximum subgraph of $G$ with no leaves, no isolated vertices, and no cycle components.
Graph $ [G] $ is called {\em the core of graph} $G$. If every component of $G$ has at most one cycle, then the core of $G$ is not defined.
 
\vskip 1.5ex
We defined a matroid $M = (E, {\cal I})$ to be {\em connected} 
 if $M$ has no loops, no coloops and if every two elements in $E$ belong to a common circuit of $M$. 

\vskip 1.5ex
We call a cocircuit $C^*$  of a connected matroid $M$ a {\em non-separating cocircuit} of $M$  if $M \setminus C^*$ is a connected matroid (see \cite{KplH}). 
 
\vskip 1.5ex

Here is some notation we will use:
\\[0.7ex]
$S (v, G)$ is the {\em $v$-star in} $G$, i.e. the set of edges incident to vertex $v$ in $G$ and $s(v, G) = |S(v, G)|$, 
\\[0.7ex]
${\mathcal S}(G) = \{S (v, G): v \in V(G)\}$ is the {\em set of vertex stars of $G$}, 
\\[0.7ex]
${\cal NC}_k^*(G)$ is the set of non-separating cocircuit of matroid $M_k(G)$, and 
\\[0.7ex]
${\cal S}_k(G)$ is the set of vertex stars $S$ of $G$ such that $|S| \le \rho_k^*(G)$, and so ${\cal S}_k(G)$ is the set of $k$-small vertex stars $S$ of $G$.

  \vskip 1.5ex 
We also remind that  graphs $G = (V, E, \phi )$ and $G' = (V', E',\phi ')$ with $E = E'$ are {\em strongly isomorphic} if
there exists a bijection $\nu :V \to V'$ such that  
$\phi (e) =\{x,y\}\Leftrightarrow \phi '(e)= \{\nu(x),\nu(y)\}$.

\vskip 1.5ex 

Matroid $M$ is called
{\em $k$-circular} 
 if there exists a graph $G$ such that M is the $k$-circular matroid of G, i.e. $M = M_k(G)$.
If $M_k(G) = M_k(G')$ implies that graphs $G$ and $G'$ are strongly isomorphic, then 
we say that $G$ is  {\em uniquely defined by $M_k(G)$}.

\vskip 1.5ex 

A natural problem $(WP)' _k$ is to describe all graphs $G$ that are uniquely defined by $M_k(G)$. In this paper we establish some properties of graphs guaranteeing that the graphs are uniquely defined by their $k$-circular matroids.

 \vskip 1ex 
  
It is easy to see that graphs $G$ and $G'$ with $E (G)= E(G')$ are  strongly isomorphic 
if and only if ${\cal S}(G) = {\cal S}(G')$. For that reason the notion of a vertex star in a graph will play a central role in 
our discussion of the problem $(WP)' _k$.

 \vskip 1ex

We will distinguish between three types of vertices according to the size of their vertex stars in a graph $G$ with respect to the rank $\rho^*_k(G)$ of matroid $M^*_k(G)$. We recall \cite{KK1} that if $M_k(G)$ is a connected matroid, then $\rho^*_k(G) = \Delta G - k + 1$.

\begin{definition}
\label{Small,Tight,Big}

{\em
Let $x$ be a vertex  in $G$, $k \ge 0$, $M_k(G)$ a connected matroid, and $\rho^*_k(G)$ the corank of $M_k(G)$. Then $x$ is called {\em a $k$-small vertex of $G$} if 
$ s(x, G) < \rho^*_k(G) + 1$, {\em a $k$-tight vertex of $G$} if 
$ s(x, G) = \rho^*_k(G) + 1$, and 
{\em a $k$-big vertex of $G$} if $ s(x, G) > \rho^*_k(G) + 1$.
}
\end{definition}

\begin{figure}[h]
\begin{center}
\scalebox{0.228}[.228]{\includegraphics{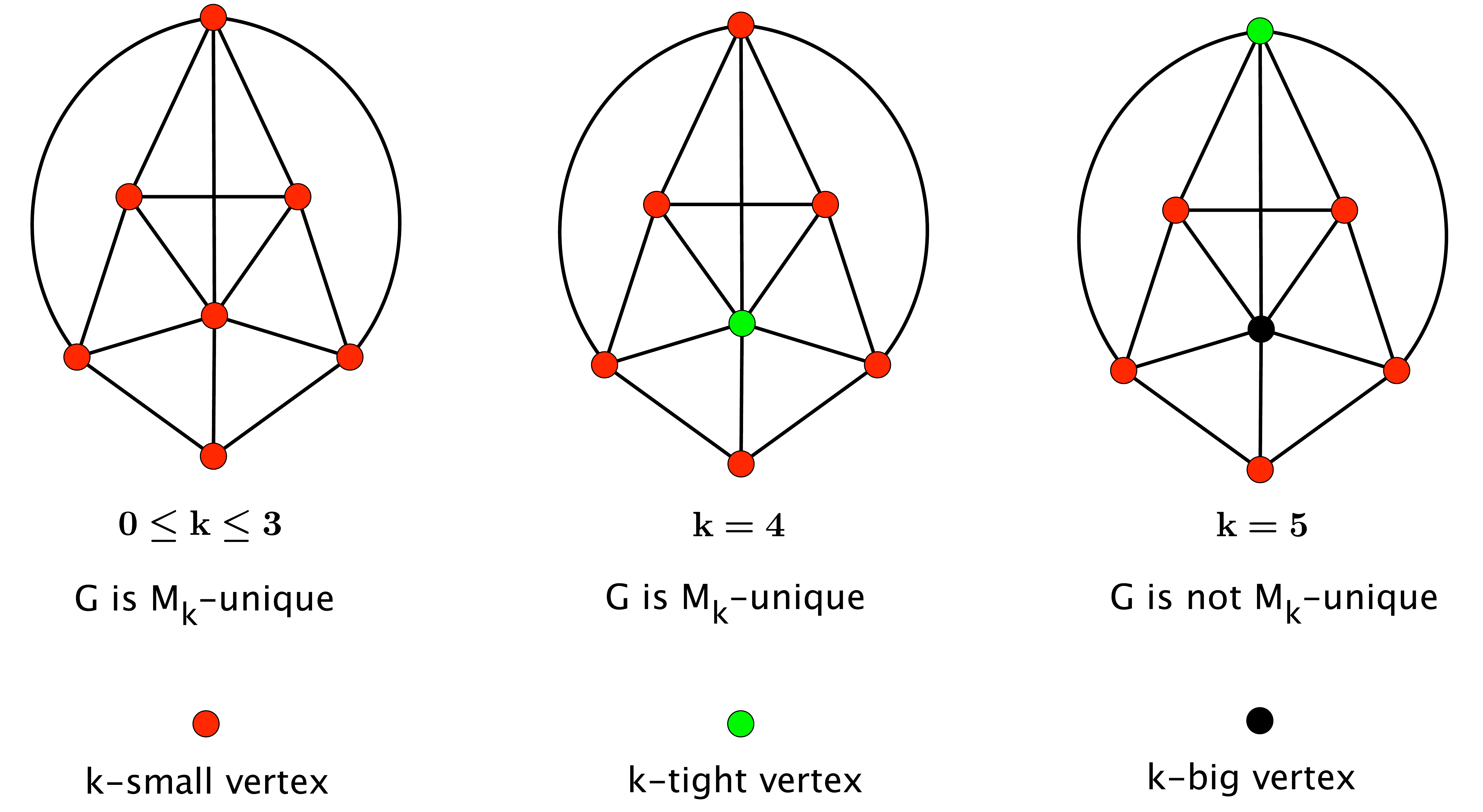}}
\end{center}
\caption
{Changes of the status of vertices when $k$ increases from 0 to 5.}
\label{k-changingG3-conFig}
\end{figure}

\begin{figure}[h]
\begin{center}
\scalebox{0.228}[.228]{\includegraphics{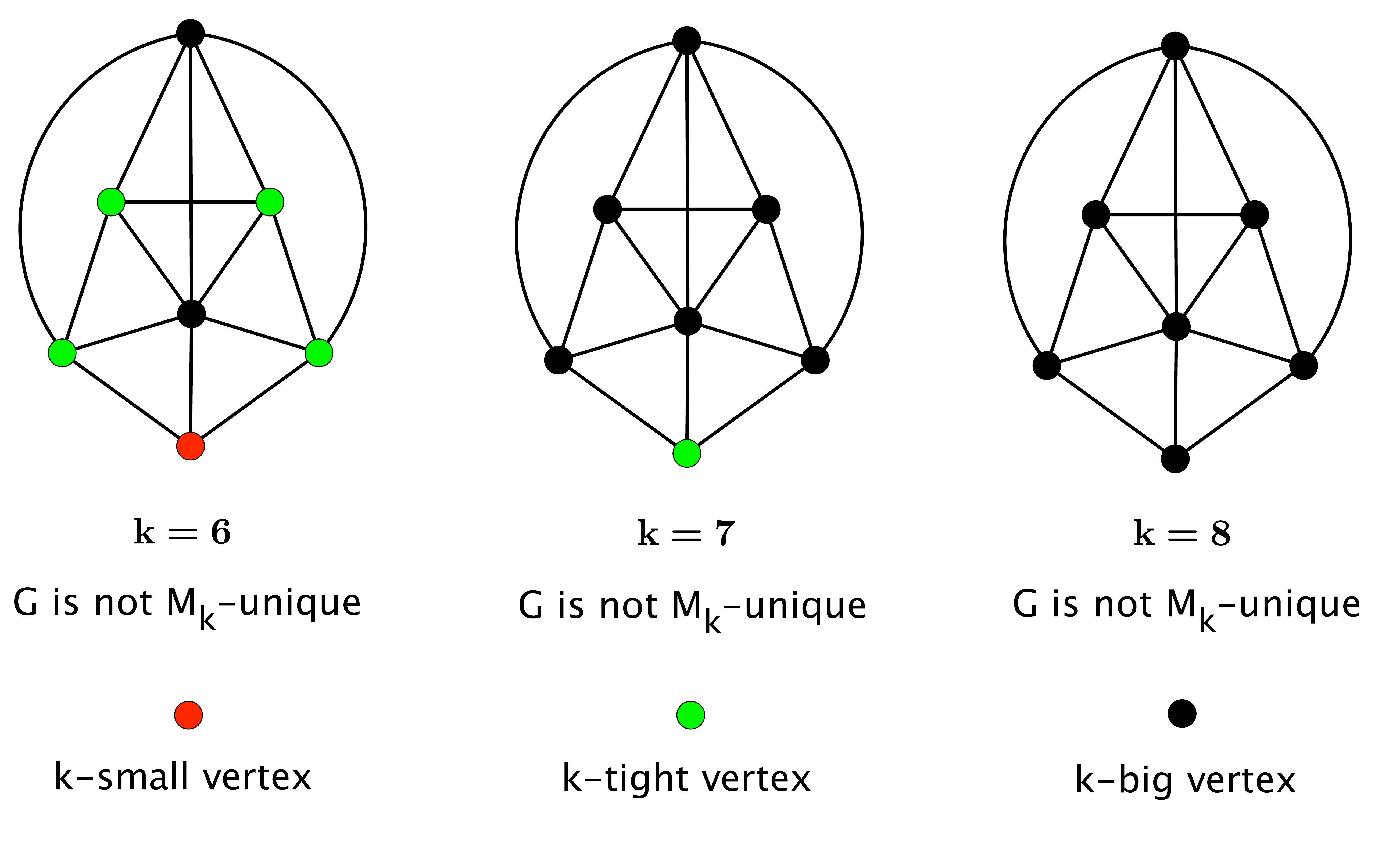}}
\end{center}
\caption
{Changes of the status of vertices when $k$ increases from 6 to 8.}
\label{k=6,7,8G3-conFig}
\end{figure}

\clearpage

\section{Stars of $G$ that are cocircuits in $M_k(G)$}
\label{Stars-Cocircuits}

\indent

We start with the following  simple and useful observation. 

\begin{claim}
\label{BigNotCo}

If $K$ is a cocircuit of a matroid, then $|K| \le \rho^*(M) + 1$.  

\end{claim}

We remind that if $B \in {\cal B}_k(G)$ and $e \in B$, then 
$K(e, B)$ denotes the fundamental cocircuit of $B$ rooted at $e$. 

\begin{claim}
\label{1StarInB}

Let $G$ be a graph, $k \ge 1$, and $x \in V(G)$. Suppose that $M_k(G)$ is a connected matroid, $B$ is a base of $M_k(G)$, $s(x, G \langle B \rangle) = 1$, and $e$ is the edge incident to 
$x$ in $G \langle B \rangle$. Then $S(x, G) = K(e, B)$, and so 
$ S(x, G) \in {\cal C}^*_k(G)$.  

\end{claim} 

\bp (uses Theorems \ref{type1} and \ref{type2})

First, suppose that $e$ is a loop in $G$. Since $s(x, G \langle B \rangle) = 1$, clearly,  $G \langle e \rangle$ is a cycle component of $G \langle B \rangle$ with exactly one vertex $x$ and one edge $e$. Then by Theorem \ref{type2}, 
 $K(e, B) \setminus e$ is the set of edges in $E \setminus B$ having $x$ as an end vertex. Therefore $S(x, G) = K(e, B)$. 

 \vskip 1ex
 
Now, suppose that $e$ is not a loop in $G$. Then $x$ is a leaf in 
$G \langle B \rangle$. By Theorem \ref{type1}, again $S(x, G) = K(e, B)$.  
\ep

\vskip 2ex

We remind  the following known fact on matroids.

\begin{claim}
\label{FundCircuit}

Let $M =  (E,{\cal I})$ be a matroid.
Then
\\[0.7ex]
$(c1)$ if $B  \in {\cal B}(M)$ and $e \in E \setminus B = B^*$, then
there exists a unique circuit $C  = C(e,B)$ of $M$ such that $e \in C \subseteq B \cup e$ $($or, equivalently, such that $C \cap B^* = \{e\}$$)$,
\\[0.7ex]
$(c2)$  similarly,
 if $B  \in {\cal B}(M)$ and $e \in B$, then
there exists a unique cocircuit $C^*  = C^*(e,B)$ of $M$ such that $e \in C^* \subseteq B^* \cup e$ $($or, equivalently, such that $C^* \cap B = \{e\}$$)$,
\\[1ex]
$(c3)$
$u \in C(e,B) \setminus e \Leftrightarrow 
(B \setminus u) \cup e \in  {\cal B}(M)$ and similarly,
$u \in C^*(e,B) \setminus e \Leftrightarrow 
(B \setminus e) \cup u \in  {\cal B}(M)$,
and
\\[1ex]
$(c4)$ for every $C \in {\cal C} (M)$ $($$C^* \in {\cal C} (M)$$)$ there exists 
$B \in {\cal B}$ and $e \in E \setminus B = B^*$ such that $C  = C(e,B)$
$($respectively, 
$e \in  B$ such that $C^*  = C^*(e,B)$$)$.

\end{claim}

\begin{claim} 
\label{LeafCocircuits}

 Let $G$ be a graph, $x$  a vertex of $G$, and $k \ge 1$. Suppose that $M_k(G)$ is a connected matroid. Then the following are equivalent:
 \\[1ex]
 $(c1)$ $S(x, G) \in {\cal C}^*_k(G)$  and  
 \\[1ex]
 $(c2)$ for every edge $e$ in $S(x, G)$ there exists $B \in {\cal B}_k(G)$ such that $e$ is either a dangling edge at $x$ or the edge of a loop component in 
$G \langle B \rangle$. 

\end{claim} 

\bp (uses Theorems \ref{type1} and \ref{type2} and Claims \ref{1StarInB} and \ref{FundCircuit})

First, we prove $(c1) \Rightarrow (c2)$. Since $S = S(x, G)$ is a cocircuit,  by Claim \ref{FundCircuit} $(c3)$, there exists $D \in {\cal B}_k(G)$ and $d \in D$ such that $S = K(d, D)$, and so 
$S \cap D = d$. Since $S$ is a vertex star of $G$, $d$ is either a dangling edge at $x$ or the edge of a loop component in $G \langle D \rangle$. Let $e \in S = K(d, D)$. Then by Theorem \ref{type1} in case edge $d$ is dangling and by Theorem \ref{type2} in case edge $d$ is the edge of a loop component in 
$G \langle D \rangle$, we  have: $B = (D \setminus d) \cup e \in {\cal B}_k(G)$. Clearly,  $K(e, B) = K(d, D) = S$. By the previous argument applied to $B$, $e$, and $S = K(e, B)$, $e$ is either a dangling edge or the edge of a loop component in $G \langle B \rangle$. 

 \vskip 1ex
 
Now, we prove $(c2) \Rightarrow (c1)$. Let $e \in  S(x, G)$ and 
$B \in {\cal B}_k(G)$ such that $e$ is either a dangling edge at $x$ or the edge of a loop component in $G \langle B \rangle$. Then $s(x, G \langle B \rangle) = 1$.  By Claim \ref{1StarInB}, $S(x, G) \in {\cal C}^*_k(G)$.
\ep

\vskip 1.5ex

We need the following three claims from \cite{KK1}.

\begin{claim}
\label{ConMGinfty}

{\em (see Claim 4.4.11 in \cite{KK1})}
Let $k \ge 1$. If $M_k(G)$ is a connected matroid, then 
$G \in {\cal G}_{\bowtie}$ and $\Delta G \ge k$.   

\end{claim}

\begin{Theorem}  {\em (4.5.3 in \cite{KK1})} {\sc Graph structure of a base of  $M_{k}(G)$ in graph $G$} 
\label{basstructure}

\vskip 0.3ex
Let $G$ be a graph and $k \geq 1$. 
Suppose that $M_k(G)$ is a connected matroid. 
Then the following are equivalent:
\\[1ex]
$(c1$) $B  \in {\cal B}_k(G)$ and
\\[1ex]
$(c2) $
$\Delta G \langle B \rangle = k - 1$, 
$V(G\langle  B \rangle) = V(G)$  $($ i.e. $B$ spans 
$V(G)$ $)$, 
and 
$\Delta A \geq 0$ 
for every component  $A$ of $ G\langle  B \rangle$
$($ i.e. $ G\langle  B \rangle$ has no tree component $)$.

\end{Theorem}

\begin{corollary}
\label{grank}
 
 {\em (4.5.4 in \cite{KK1})}
Let $G$ be a graph and $k \geq 1$. Suppose that $M_k(G)$ is a connected matroid. 
Then 
$\rho _k(G) = |V(G)|  - 1 + k $ and  
 $\rho ^*_k(G) = |E(G)| - |V(G)| + 1 - k$.

\end{corollary}

\begin{claim}
\label{SmallTight}

Let $G$ be a graph and $k \ge 1$. Suppose that 
\\[1ex]
 $(a1)$ $M_k(G)$ is a connected matroid and
\\[1ex]
 $(a2)$  $x$ is a vertex of $G$ such that $ s(x, G) \le \rho^*_k(G) + 1$, i.e. $x$ is not a $k$-big vertex of $G$.

 \vskip 1ex

Then the following are equivalent:
\\[1ex]
$(c1)$ there exists $B \in {\cal B}_k(G)$ such that 
$x \notin V \lfloor G \langle B \rangle \rfloor  $ and 
\\[1ex]
$(c2)$ $Q \setminus x$ has a cycle, where $Q$ is the component of $G$ containing $x$. 

\end{claim} 

\bp (uses Theorem \ref{basstructure}, Corollary \ref{grank}, and Claim \ref{ConMGinfty}) 
\\[1.5ex]
${\bf (p1)}$ We prove $(c1) \Rightarrow (c2)$. Suppose that $B \in {\cal B}_k(G)$ and $x \notin V  \lfloor G \langle B \rangle \rfloor$. By Theorem \ref{basstructure}, every component of $ G \langle B \rangle $ has a cycle. Let  $A$ be the component of $G \langle B \rangle$ containing $x$. Since 
$x \notin V \lfloor G \langle B \rangle \rfloor$, clearly, $A$ has a cycle $C$ avoiding $x$. Obviously, $Q$ contains $A$, and therefore also contains $C$.  
\\[1ex]
${\bf (p2)}$ We prove $(c2) \Rightarrow (c1)$. 
Since $M_k(G)$ is connected, by Claim \ref{ConMGinfty}, $G \in {\cal G}_{\bowtie}$. Since $G \in {\cal G}_{\bowtie}$, clearly, $G \setminus e$ has no tree component for every edge $e$ in $G$. Then there exists a maximal subset $Z$ of $S = S(x, G)$ such that $G' = G \setminus Z$ has no tree  component and there is an edge $p$ in $S \setminus Z$ incident to a component of $G \setminus x$ having a cycle. By assumption $(c2)$, such edge set $Z$ exists and by the maximality of $Z$ edge $p$ is unique. Therefore $a \notin E \lfloor G' \rfloor$ for every $a \in S \setminus Z  = S \cap E(G')$. Thus, $x \notin V \lfloor A' \rfloor$ for every $A' \in Cmp( G')$.

 \vskip 1ex
 
Since $G \setminus e$ has no tree component for every edge $e$ in $G$, set 
$Z$ is a proper subset of $S$. Then, clearly, 
$V(G') = V(G)$. Since $M_k(G)$ is a connected matroid, by Corollary \ref{grank}, we have: $\rho ^*_k(G) = |E(G)| - |V(G)| + 1 - k$. Since 
$ s(x, G) \le \rho^*_k(G) + 1$, we have: 
$ |E(G)|  - s(x, G) + 1 - |V(G)| \ge k - 1$. Since $|Z| <  s(x, G)$, we have: 
$\Delta G' = |E(G')| - |V(G')| = |E(G)| - |Z| - |V(G)| \ge |E(G)|  - s(x, G) + 1 - |V(G)| \ge k - 1$. 

 \vskip 1ex
 
Let ${\cal F}$ be the family of subgraphs $F$ of $G'$ such that $V(F) = V (G')$, $\Delta F \ge k-1$, and $F$ has no tree component. Let 
$H \in {\cal F}$ be such that $\Delta H$ is minimum. Let $B = E(H)$. Since $H$ is a subgraph of $G'$, clearly $x \notin V \lfloor A \rfloor$ for every 
$A \in Cmp( H ) = Cmp( G \langle B \rangle )$. We claim that $B \in {\cal B}_k(G)$. By Theorem \ref{basstructure}, it is sufficient to prove that $\Delta H = k-1$. 

 \vskip 1ex
 
First, suppose that every component of $H$ has exactly one cycle. Then 
$\Delta H = 0$. Since $0 = \Delta H \ge k-1$ and $k$ is a postive integer, clearly $k = 1$, and therefore $\Delta H = k-1$. 

 \vskip 1ex
 
Now, suppose that $H$ has a component $D$ such that $D$ has at least two cycles. Let $e$ be an edge that belongs to a cycle in $D$. Let $H' = H \setminus e$. Then, clearly, $V(H') = V (G')$, $H'$ has no tree component, and 
$\Delta H' = \Delta H - 1$. By the selection of $H$, graph $H' \notin {\cal F}$, and so $\Delta H' < k-1$. Therefore, $\Delta H = k-1$. 
\ep

\begin{claim}
\label{A-e}

{\em (3.3.9 in \cite{KK1})}
Let $A$ be a connected graph with at least one cycle and 
$e \in E(A)$.

\vskip 0.3ex
\noindent
Then
\\[1ex]
$(c0)$ $e \in E \lfloor A \rfloor $ if and only if both end vertices of $e$ belong to $ \lfloor A \rfloor $, 
\\[1ex]
$(c1)$  if $e \notin E \lfloor A \rfloor $, then $A \setminus e$ has  two components and exactly one of them is a tree and the other component contains $\lfloor A \rfloor $,
\\[1ex]
$(c2)$  if $A$ has one cycle and $e \in E \lfloor A \rfloor $, 
then $A \setminus e$  is a tree,
\\[1ex]
$(c3)$ if $A$ has at least two cycles and  $e \in E \lfloor A \rfloor $, then  
every component of $A \setminus e$ contains a cycle, and
\\[1ex]
$(c4)$ if $A \setminus e$ has two components and 
$v \in V \lfloor A \rfloor $, then the component of $A \setminus e$ containing $v$ is not a tree.

\end{claim}

\begin{notation}
\label{GvDefin}

{\em 
Let $G$ be a graph and $x \in V(G)$. If $G \setminus x$ has a cycle, then let $G_x$ denote the union of all non-tree components of $G \setminus x$. If 
$G \setminus x$ has a tree component, then let $G^x$ denote the union of all tree components of $G \setminus x$. 
}

\end{notation}

\begin{claim} 
\label{OutFloorGeneral}

Let $G$ be a graph, $x \in V(G)$, and $k \ge 1$. Suppose that 
\\[1ex]
$(a1)$ $M_k(G)$ is a connected matroid and   
\\[1ex]
$(a2)$ there exists $B \in {\cal B}_k(G)$ such that 
$x \notin V \lfloor G \langle B \rangle \rfloor $. 

\vskip 1ex
Then there exists $B' \in {\cal B}_k(G)$ such that 
\\[1ex]
$(c1)$ $T$ is a tree component of 
$G \langle B' \rangle \setminus x$ if and only if $T$ is a tree component of 
$G \setminus x$, 
\\[1ex]
$(c2)$ $G \langle B' \rangle \setminus V(G_x)$ is a tree, and 
\\[1ex]
$(c3)$ there is exactly one edge in $G \langle B' \rangle$ from $x$ to $V(G_x)$.    

\end{claim} 

\bp (uses Theorem \ref{basstructure} and \ref{type1}  and Claims \ref{A-e} and \ref{FundCircuit})

Let 
${\cal F} = \{ B \in {\cal B}_k(G): x \notin V \lfloor G \langle B \rangle \rfloor\}$. By assumption $(a2)$ of our claim, ${\cal F} \ne \emptyset$.  Let 
$B' \in {\cal F}$ be such that $s(x,  G \langle B' \rangle)$ is minimum. Since $B' \in {\cal B}_k(G)$, by Theorem \ref{basstructure}, $B'$ spans $V(G)$ and $G \langle B' \rangle$ has no tree components. We claim that 
$G \langle B' \rangle$ satisfies conditions 
$(c1)-(c3)$ of our claim.

Since  by $(a2)$, $x \notin V \lfloor G \langle B' \rangle \rfloor$, graph $G  \setminus x$ has a cycle and $G_x$ is defined. 

 \vskip 1ex

Let $A'$ be the component of $G \langle B' \rangle $ containing $x$, and so 
$x \notin V \lfloor A' \rfloor$. Since $x$ is not incident to  $\lfloor A' \rfloor$, clearly,  $ e \notin E \lfloor A' \rfloor $ for every $e \in S(x, B')$. 

By Claim \ref{A-e} $(c2)$, for every $e \in S(x, A')$ the graph
 $A' \setminus e$ consists of two components: 
one of these two components  is a tree $T_e$ and the other component 
$U_e$ contains $\lfloor A' \rfloor$. 

Let $P$ be a minimal path in $A'$ from $x$ to $\lfloor A' \rfloor$. We claim that 
$P$ is the only minimal path in $A'$ from $x$ to $\lfloor A' \rfloor$. Indeed, if there is another minimal path $P'$ from $x$ to $\lfloor A' \rfloor$, then the core of 
$P \cup P'  \cup \lfloor A' \rfloor$ is a subgraph of $A'$ containing 
$\lfloor A ''\rfloor$ properly, contradicting the definition of $\lfloor A' \rfloor$. Obviously, $S(x,A')$ and $E(P)$ have exactly one edge in common, say $f$.    

Suppose that there exists a tree component $T$ of $G \langle B' \rangle \setminus x$. Then clearly, $T =T_e$ for some $e \in S(x, A') \setminus f$. We claim that $T$ is a tree component of $G \setminus x$. Indeed, suppose not. Then there exists $g \in E\setminus B'$ such that $g$ is incident to $T$ and not incident to $x$. Therefore by Theorem \ref{type1}, $g \in K(e, B')$. By Claim \ref{FundCircuit},  $B'' = (B' \setminus e) \cup g$ is a $k$-base of $G$. Clearly, $x \notin V \lfloor G \langle B'' \rangle \rfloor$, and so 
$B'' \in {\cal F}$. But $ s(x, G \langle B'' \rangle ) <  s(x, G \langle B' \rangle)$, contradicting the selection of $B'$.  

Now, suppose that there exists a tree component $T$ of $G \setminus x$. 
Then every cycle of $G$ having vertices of $T$ contains $x$. Since $B'$ spans $V(G)$ and  
$G \langle B' \rangle$ has no tree components, there exists an edge 
$e \in B'$ incident to a vertex $z$ of $T$ and $V(G) \setminus V(T)$. Since $T$ is a tree component of $G \setminus x$, edge $e$ is incident to vertex 
$x$. Since $T_e$ is a component of $G \langle B' \rangle \setminus e$, $T_e$ is a subgraph of the component of $G \setminus x$ containing $z$. Therefore $T_e$ is a subtree of $T$. By the argument in the previous paragraph $T_e$ is a tree component of $G \setminus x$. Therefore  
$T_e = T$ and $e$ is the only edge in $B'$ having exactly one end vertex in $T$. Thus, our claims $(c1)$ and $(c2)$ are true. 

Since $T$ is a tree component of $G \langle B' \rangle \setminus x$ if and only if $T$ is a tree component of $G \setminus x$, every vertex of $V(G_x)$ belongs to a non-tree component of $G \langle B' \rangle \setminus x$. If there are two edges in $G \langle B' \rangle$ from $x$ to $V(G_x)$, then $x \in [ G \langle B' \rangle ]$, contradicting that 
$x \notin V \lfloor G \langle B' \rangle \rfloor $. Therefore $f$ is the only edge in $G \langle B' \rangle$ from $x$ to $V(G_x)$, and so $(c3)$ of our claim is true.  
\ep

\vskip 3.5ex

\begin{claim} 
\label{OutFloorToLeaf}

Let $G$ be a graph and $k \ge 1$. Suppose that 
\\[1ex]
$(a1)$ $M_k(G)$ is a connected matroid, 
\\[1ex]
$(a2)$ there exists $B \in {\cal B}_k(G)$ and $x \in V(G)$ such that 
$x \notin V \lfloor G \langle B \rangle \rfloor $, and  
\\[1ex]
$(a3)$ $G \setminus x$ has no tree component. 

\vskip 1ex

Then there exists $B' \in {\cal B}_k(G)$ such that $x$ is a leaf in 
$G \langle B' \rangle$.   

\end{claim} 

\bp (uses Claims \ref{OutFloorGeneral})

The assumptions $(a1)$ and $(a2)$ of our claim are the assumptions $(a1)$ and $(a2)$, respectively of Claim \ref{OutFloorGeneral}. Therefore there exists $B' \in {\cal B}_k(G)$ satisfying claims $(c1)-(c3)$ of Claim \ref{OutFloorGeneral}. Since $G \setminus x$ has no tree component, by Claim \ref{OutFloorGeneral} $(c1)$, graph $G \langle B' \rangle \setminus x$ has no tree component. Therefore $V(G) = V(G_x) \cup \{x\}$. By Claim \ref{OutFloorGeneral} $(c2)$, $G \langle B' \rangle \setminus V(G_x)$ is the isolated vertex $x$. By Claim \ref{OutFloorGeneral} $(c3)$, there is exactly one edge in $G \langle B' \rangle$ from $x$ to $V(G_x)$, and so $x$ is a leaf in $G \langle B' \rangle$.  
\ep

\vskip 1.5ex

From Claims \ref{SmallTight} and \ref{OutFloorToLeaf}, we have:

\begin{claim}
\label{SmallTightLeaf}

Let $G$ be a graph and $k \ge 1$. Suppose that 
\\[1ex]
$(a1)$ $M_k(G)$ is a connected matorid, 
\\[1ex]
$(a2)$ $x$ is a vertex of $G$ such that $ s(x, G) \le \rho^*_k(G) + 1$, i.e. $x$ is not a $k$-big vertex of $G$, and
\\[1ex]
$(a3)$ $G \setminus x$ has no tree component.  

\vskip 1ex 

Then exactly one of the following holds:
\\[1ex]
$(c1)$ the component of $G$ containing vertex $x$ has no other vertex, 
\\[1ex]
$(c2)$ there exists $B \in {\cal B}_k(G)$ such that $x$ is a leaf in 
           $G \langle B \rangle$.

\end{claim}

From Claims \ref{LeafCocircuits} and \ref{SmallTightLeaf}, we have: 

\begin{Theorem} {\sc A condition for a star to be a $k$-cocircuit}
\label{SmallTightCo}

Let $G$ be a graph, $k \ge 1$, and $x \in V(G)$. Suppose that 
\\[1ex]
$(a1)$ $M_k(G)$ is a connected matroid and
\\[1ex]
$(a2)$ $ s(x, G) \le \rho^*_k(G) + 1$, i.e. $x$ is not a $k$-big vertex of $G$.

\vskip 1ex

Then the following are equivalent:
\\[1ex]
$(c1)$ $S(x, G) \in {\cal C}^*_k(G)$ and
\\[1ex]
$(c2)$ $G \setminus x$ has no tree component. 

\end{Theorem}

\bp (uses Theorem \ref{basstructure}, Corollary \ref{grank}, and Claims \ref{LeafCocircuits}, \ref{SmallTight}, and \ref{SmallTightLeaf}) 
\\[1ex]
{\bf (p1)}
First, we prove $(c1) \Rightarrow (c2)$. Let $e \in S(x, G)$. Since 
$S(x, G) \in {\cal C}^*_k(G)$, by Claim \ref{LeafCocircuits}, there exists 
$B \in {\cal B}_k(G)$ such that $e$ is either a dangling edge or the edge of a loop component in $G \langle B \rangle$. By Theorem \ref{basstructure}, every component in $G \langle B \rangle$ contains a cycle. Therefore every component of $G \langle B \rangle \setminus x$ has a cycle. Thus, $G \setminus x$ has no tree component. 
\\[1ex]
{\bf (p2)} 
Now, we prove $(c2) \Rightarrow (c1)$. If the component of $G$ containing $x$ has at least two vertices, by Claim \ref{SmallTightLeaf}, there exists 
$B \in {\cal B}_k(G)$ such that $x$ is a leaf in $G \langle B \rangle$. Therefore, by Claim \ref{LeafCocircuits}, $S = S(x, G) \in {\cal C}^*_k(G)$. So we assume that the component of $G$ containing vertex $x$ has no other vertex, and so every edge in $S(x, G)$ is a loop in $G$.

In what follows, the arguments we use are similar to those in ${\bf (p2)}$ of the proof of Claim \ref{SmallTight}.

Let $e \in S = S(x, G)$ and $ Z = S \setminus e$. Let 
$G' = G \setminus Z$. Since $M_k(G)$ is a connected matroid, by Corollary \ref{grank}, we have: $\rho ^*_k(G) = |E(G)| - |V(G)| + 1 - k$. Since 
$ s(x, G) \le \rho^*_k(G) + 1$, we have: $ |E(G)|  - s(x, G) + 1 - |V(G)| \ge k - 1$. Since $|Z| <  s(x, G)$, we have: 
$\Delta G' = |E(G')| - |V(G')| = |E(G)| - |Z| - |V(G)| \ge |E(G)|  - s(x, G) + 1 - |V(G)| \ge k - 1$. Clearly $V (G') = V(G)$, $\Delta G' \ge k-1$, and $G'$ has no tree component.

\vskip 1ex

Let ${\cal F}$ be the family of subgraphs $F$ of $G'$ such that $V(F) = V (G')$, $\Delta F \ge k-1$, and $F$ has no tree component. Let 
$H \in {\cal F}$ be such that $\Delta H$ is minimum. Let $B = E(H)$. We claim that $B \in {\cal B}_k(G)$. By Theorem \ref{basstructure}, it is sufficient to prove that $\Delta H = k-1$. 

\vskip 1ex

First, suppose that every component of $H$ has exactly one cycle. Then 
$\Delta H = 0$. Since $0 = \Delta H \ge k-1$ and $k$ is a postive integer, clearly $k = 1$, and therefore $\Delta H = k-1$.

Now, suppose that $H$ has a component $D$ such that $D$ has at least two cycles. Let $p$ be an edge that belongs to a cycle in $D$. Let 
$H' = H \setminus p$. Then clearly, $V(H') = V (G')$, $H'$ has no tree component, and $\Delta H' = \Delta H - 1$. By the selection of $H$, graph 
$H' \notin {\cal F}$, and so $\Delta H' < k-1$. Therefore, $\Delta H = k-1$. 

\vskip 1ex

Thus, in both cases $\Delta H = k-1$, and so $B \in {\cal B}_k(G)$. By definition of $H$ and since the component of $G$ containing vertex $x$ has no other vertex, clearly $e$ is the only edge of the component of $H = G \langle B \rangle$ containing $x$. By Claim \ref{LeafCocircuits}, $S(x, G) \in {\cal C}^*_k(G)$.  
\ep

\section{Stars of $G$ and non-separating cocircuits of $M_k(G)$}
\label{NonSep}

We will use the following three facts from \cite{KK1}. 

\begin{Theorem} {\em (4.2.1 in \cite{KK1})} {\sc  A criterion for matroid $M_{k}(G)$ to be non-trivial}
\label{maxk}

\vskip 0.3ex

Let $G$ be a graph, $F(G)$ the union of all tree components of $G$,
 and $k \geq 1$. Then the following are equivalent:    
\\
$(a1)$ ${\cal C}_k(G) \ne \emptyset$,       
\\
$(a2)$ $k \leq \Delta G + cmp(F(G))$, and 
\\
$(a3)$ $M_k(G)$ is a non-trivial matroid.     

\end{Theorem}

\begin{Theorem} {\em (4.4.5 in \cite{KK1})} 
{\sc Graph description of connected matroid $M_k(G)$}
\label{MkCon}

\vskip 0.3ex
 
Let $k \ge 2$. Then the following are equivalent:
\\[1ex]
$(a1)$ $M_k(G)$ is a non-trivial matroid and $G \in {\cal G}_{\bowtie}$ and 
\\[1ex]
$(a2)$ $M_k(G)$ is a connected matroid. 

\end{Theorem}

\begin{Theorem} {\em (4.4.6 in \cite{KK1})} 
{\sc Graph description of connected matroid $M_1(G)$}
\label{M1Con}

\vskip 0.3ex
The following are equivalent:
\\[1ex]
$(a1)$ $M_1(G)$ is a non-trivial matroid and $G \in {\cal CG}_{\bowtie}$ and 
\\[1ex]
$(a2)$ $M_1(G)$ is a connected matroid. 

\end{Theorem}

Recall that a cocircuit $K$ of a connected matroid $M$ is {\em non-separating} if and only if matroid $M \setminus K$ is connected. Given a graph $G$, let ${\cal NC}_k^*(G)$ denote the set of non-separating cocircuit of $M_k(G)$.

\vskip 1.5ex

Here is a graph description of non-separating cocircuits of $M_k(G)$ in terms of $k$ and  $G$.

\begin{claim}
\label{k-NonSepInG}

Let $G$ be a graph, $K \subseteq  E(G)$, and $k \ge 1$.
Then the following are equivalent:
\\[1ex]
$(c1)$  $K \in {\cal NC}_k^*(G)$ and 
\\[1ex]
$(c2)$ $\Delta (G \setminus K) \ge k$ and 

\vskip 0.7ex
$(c2.1)$ $G \setminus K \in {\cal CG}_{\bowtie}$ if $k = 1$ and 
\vskip 0.7ex
$(c2.2)$ $G \setminus K \in {\cal G}_{\bowtie}$ if $k \ge 2 $.

\end{claim}

\bp (uses Theorem \ref{maxk}, \ref{MkCon}, and \ref{M1Con})
 
By definition, $K \in {\cal NC}_k^*(G)$ if and only if $M_k(G) \setminus K$ is a connected matroid. Clearly, $M_k(G) \setminus K = M_k(G \setminus K)$. Therefore $K \in {\cal NC}_k^*(G)$ if and only if matroid $M_k(G \setminus K)$ is connected.

\vskip 1ex

Suppose that $k = 1$. By Theorem \ref{M1Con} for $k = 1$, 
$M_1(G \setminus K)$ is connected if and only if $M_1(G)$ is non-trivial  and 
$G \in {\cal CG}_{\bowtie}$. 

\vskip 1ex

Next, suppose that $k \ge 2$. By Theorem \ref{MkCon} for $k \ge 2$, 
$M_k(G \setminus K)$ is connected if and only if $M_k(G)$ is non-trivial  and 
$G \in {\cal CG}_{\bowtie}$.

\vskip 1ex

Now, if $G \setminus K \in {\cal G}_{\bowtie}$, then $G \setminus K$ has no tree components and by Theorem \ref{maxk}, $M_k(G)$ is non-trivial if and only if $\Delta (G \setminus K) \ge k$. Thus, $(c1) \Leftrightarrow (c2)$.
 \ep

\vskip 2ex

We will distinguish between three possible types of $(B,e)$-cocircuits $K(e,B)$ depending on the structure of component $A$ of $ G \langle B \rangle $  containing edge $e$ and on the position of edge $e$ in $A$.
\begin{definition} {\em$($4.7.1 in \cite{KK1}$)$
\label{cocircuits-type}

\vskip 0.3ex
Let $B \in {\cal B}_k(G)$, $e \in B$, and $k \ge 1$. Then
\\[0.7ex]
$(t1)$ $K(e,B)$ is a {\em $(B,e)$-cocircuit  in $M_k(G)$ of type 1},
if $e \notin E \lfloor A \rfloor $, where $A$ is a component of 
$ G \langle B \rangle$ containing edge $e$,
\\[0.7ex]
$(t2)$ $K(e,B)$ is a {\em $(B,e)$-cocircuit  in $M_k(G)$ of type 2},
if $e \in E \lfloor A \rfloor $, where $A$ is a unicyclic component of 
$ G \langle B \rangle$ containing edge $e$, and
 \\[0.7ex]
$(t3)$ $K(e,B)$ is a {\em $(B,e)$-cocircuit  in $M_k(G)$ of type 3},
if $e \in E [ A ] $, where $A$ is a component of 
$ G \langle B \rangle$ that has at least two cycles and contains edge $e$.}
 \end{definition}

 \begin{Theorem} {\em(4.7.2  in \cite{KK1})} {\sc Graph description of rooted cocircuits of type 1}
 \label{type1}
 
\vskip 0.3ex
Let $G$ be a graph and $k \ge 1$. Suppose that  
\\[1ex] 
$(a1)$ $M_k(G)$ is a connected matroid, $B \in {\cal B}_k(G)$, and $e \in B$ and  
\\[1ex] 
$(a2)$  $K(e,B)$ is a $(B,e)$-cocircuit  in $M_k(G)$ of type 1 $($i.e. edge $e \notin E \lfloor A \rfloor $, where $A$ is a component of 
$ G \langle B \rangle$ containing edge $e$$)$.  

\vskip 1ex
Then exactly one of the two components of $A \setminus e$ is a tree $T$ 
and $K(e,B) = K'(e,B) \cup e$, where $K'(e,B)$ is the set of edges in  
$E \setminus B$ having at least one end-vertex in $V(T)$.

 \end{Theorem}

\begin{Theorem} {\em(4.7.3 in \cite{KK1})} {\sc Graph description of rooted cocircuits of type 2}
 \label{type2}
 
\vskip 0.3ex
Let $G$ be a graph and $k \ge 1$. Suppose that  
\\[1ex] 
$(a1)$ $M_k(G)$ is a connected matroid, $B \in {\cal B}_k(G)$, and $e \in B$ and  
\\[1ex] 
$(a2)$ $K(e,B)$ is a $(B,e)$-cocircuit  in $M_k(G)$ of type 2, $($i.e. 
$e \in E \lfloor A \rfloor $, where $A$ is a unicyclic component of 
$ G \langle B \rangle$ containing edge $e$$)$.    

\vskip 1ex
Then $A \setminus e$ is a tree and $K(e,B) = K'(e,B) \cup e$, where $K'(e,B)$ is the set of edges in $E \setminus B$ having at least one end-vertex in 
$V(A \setminus e)$.

 \end{Theorem}

\vskip 1ex
\begin{Theorem}{\em(4.7.4 in \cite{KK1})} {\sc Graph description of rooted cocircuits of type 3}
 \label{type3}
 
\vskip 0.3ex
Let $G$ be a graph and $k \ge 1$. Suppose that  
\\[1ex] 
$(a1)$ $M_k(G)$ is a connected matroid, $B \in {\cal B}_k(G)$, and $e \in B$ and  
\\[1ex] 
$(a2)$ $K(e,B)$ is a $(B,e)$-cocircuit  in $M_k(G)$ of type 3, $($i.e. 
$e \in E [ A ] $, where $A$ is a component of $ G \langle B \rangle$ that has at least two cycles and contains edge $e$$)$.    

\vskip 1ex
Then $K(e,B) = (E \setminus B) \cup e$.  
 
 \end{Theorem}

\begin{claim}
\label{NonSepToStar}

Let $G$ be a graph and $k \ge 1$. If $M_k(G)$ is a connected matroid, then 
\vskip 0.3ex
\noindent
${\cal NC}_k^*(G) \subseteq S(G)$.     

\end{claim} 

\bp (uses Theorems \ref{type1}, \ref{type2}, and \ref{type3} and Claim \ref{FundCircuit})

\vskip 0.3ex
Since $K \in {\cal C}_k^*(G)$, by Claim \ref{FundCircuit} $(c3)$,  there exists 
$B \in {\cal B}_k(G)$ and $e \in B$ such that $K$ is the $(B, e)$-cocircuit in $M_k(G)$.  
 
 \vskip 1ex
 
 First, suppose that  $K$ is the $(B, e)$-cocircuit in $M_k(G)$ of type 3.
Then  by Definition \ref{cocircuits-type}, $e \notin E \lfloor A \rfloor $, where $A$ is a component of 
$ G \langle B \rangle$ containing edge $e$. Hence by Theorem \ref{type3}, $K = (E(G) \setminus B) \cup e$, and so $E(G) \setminus K \subset B$. Therefore $E(G) \setminus K \in {\cal I}_k(G)$. Thus, $M_k(G) \setminus K$ is not a connected matroid, contradicting  $K \in {\cal NC}_k^*(G)$.  

 \vskip 1ex
 
Next,  suppose that  $K$ is the $(B, e)$-cocircuit in $M_k(G)$ of type 2.  
%I 
Then by Definition \ref{cocircuits-type}, $e \in E \lfloor A \rfloor $, where $A$ is a unicyclic component of $ G \langle B \rangle$ containing edge $e$. Hence by Theorem \ref{type2}, $A \setminus e$ is a tree and $K \setminus e$ is the set of edges in $E \setminus B$ having at least one end-vertex in $V(A \setminus e)$. If $A$ has at least two edges, then $A \setminus e$ is a tree component of $G \setminus K$
having at least one edge. Clearly, no edge of $A \setminus e$ belongs to a $k$-circuit of $G \setminus K$. Therefore  $M_k(G \setminus K)$ is not a connected matroid, again contradicting  $K \in {\cal NC}_k^*(G)$. Therefore $e$ is the only edge of the component of $G \langle B \rangle$ and it is a loop incident to exactly one vertex, say $x$. Clearly, $S(x, G)$ is the set of edges in $E \setminus B$ incident to $x$, and so $K= S(x, G)$. 

 \vskip 1ex
 
Finally,   suppose that  $K$ is the $(B, e)$-cocircuit in $M_k(G)$ of type 1.       
Then by Definition \ref{cocircuits-type}, $e \in E [ A ] $, where $A$ is a component of 
$ G \langle B \rangle$ that has at least two cycles and contains edge $e$. Hence by Theorem \ref{type1}, exactly one of the two components of $A \setminus e$ is a tree $T$ and $K \setminus e$ is the set of edges in  $E \setminus B$ having at least one end-vertex in $V(T)$. Therefore $T$ is a tree component of $G \setminus K$. If $|V(T)| \ge 2 $, then $|E(T)| \ge 1$ and no edge of $T$ belongs to a $k$-circuit of $G \setminus K$. Therefore  $M_k(G \setminus K)$ is not a connected matroid, contradicting  $K \in {\cal NC}_k^*(G)$. Thus, $T$ has exactly one vertex, say $x$. Clearly, is $S(x, G)$ is the set of edges in $E \setminus B$ incident to $x$, and so $K= S(x, G)$.     

 \vskip 1ex
 
Thus, we have proved that if $K \in {\cal NC}_k^*(G)$, then there exists $x \in V(G)$ such that $K= S(x, G)$. 
\ep

\begin{claim}
\label{SxSmallMG-xNonTriv}
 
Let $G$ be a graph, $x$ a vertex of $G$, and $k \ge 1$. Suppose that
\\[1ex]
$(a1)$ $M_k(G)$ is a connected matroid and 
\\[1ex] 
$(a2)$ $s(x, G) \le \rho^*_k(G)$, i.e. $x$ is a $k$-small vertex of $G$. 

\vskip 1ex 
 Then $M_k(G \setminus x)$ is a non-trivial matroid.   

\end{claim} 

 \bp (uses Theorem \ref{maxk} and Corollary \ref{grank})  
 
 First, we prove that $M_k(G \setminus x)$ is a non-trivial matroid. Since $M_k(G)$ is a connected matroid, by Corollary \ref{grank}, we have: 
$\rho ^*_k(G) = |E(G)| - |V(G)| + 1 - k$. Since $ s(x, G) \le \rho^*_k(G)$, we have: 
$ k  \le (|E(G)|  - s(x, G)) - (|V(G)| - 1) = |E(G \setminus x)| - |V(G \setminus x)| = \Delta (G \setminus x) $. 
Therefore $\Delta (G \setminus x) \ge k$. Therefore by Theorem \ref{maxk}, 
$M_k(G \setminus x)$ is a non-trivial matroid. 
\ep

\begin{claim} 
\label{G-xInftyNonSep}

Let $G$ be a graph, $x$ a vertex of $G$, and $k \ge 1$. Suppose that
\\[1ex]
$(a1)$ $M_k(G)$ is a connected matroid,
\\[1ex]
$(a2)$ $s(x, G) \le \rho^*_k(G)$, i.e. $x$ is a $k$-small vertex of $G$, and  
\\[1ex] 
$(a3)$ $G \setminus x \in {\cal CG}_{\bowtie}$ if $k = 1$ and 
$G \setminus x \in {\cal G}_{\bowtie}$ if $k \ge 2$. 

\vskip 1ex 
 Then $S(x, G) \in {\cal NC}_k^*(G)$.    

\end{claim}

\bp (uses Theorems \ref{SmallTightCo}, \ref{MkCon}, and  \ref{M1Con} and Claim \ref{SxSmallMG-xNonTriv}) 

By assumption $(a3)$, $G \setminus x$ has no tree components. Therefore by Theorem \ref{SmallTightCo}, 
\\
$S = S(x, G) \in {\cal C}_k(G)$. Clearly, 
$M_k(G) \setminus S = M_k(G \setminus x)$. We need to prove that $M_k(G \setminus x)$ is a connected matroid. By Claim \ref{SxSmallMG-xNonTriv}, $M_k(G \setminus x)$ is a non-trivial matroid. 

 \vskip 1ex
 
Suppose that $k = 1$. Then since $G \setminus x \in {\cal CG}_{\bowtie}$, by Theorem \ref{M1Con}, $M_k(G \setminus x)$ is a connected matroid, and so 
$S(x, G) \in {\cal NC}_k^*(G)$.

 \vskip 1ex
 
Now, suppose that $k \ge 2$. Then since $G \setminus x \in {\cal G}_{\bowtie}$, by Theorem \ref{MkCon}, $M_k(G \setminus x)$ is a connected matroid, and so 
$S(x, G) \in {\cal NC}_k^*(G)$. 
\ep

 \vskip 2ex

\begin{Theorem} {\sc Graph structure of non-separating $k$-cocircuits} 
\label{k>1NonSepStructure}
 
 \vskip 0.3ex
Let $G$ be a graph and $k \ge 1$. Suppose that $M_k(G)$ is a connected matroid. Then the following are equivalent: 
\\[1ex] 
$(c1)$ $K \in {\cal NC}_k^*(G)$ and 
\\[1ex] 
$(c2)$ there exists $x \in V(G)$ such that $K= S(x, G)$,  
$s(x, G) \le \rho^*_k(G)$, and 

\vskip 0.6ex
$(c2.1)$ $G \setminus x \in {\cal CG}_{\bowtie}$ if $k = 1$ and 
\vskip 0.6ex
$(c2.2)$ $G \setminus x \in {\cal G}_{\bowtie}$ if $k \ge 2$. 

\end{Theorem} 

\bp (uses Theorem \ref{MkCon} and \ref{M1Con} and Claims \ref{NonSepToStar} and \ref{G-xInftyNonSep}) 

By Claim \ref{G-xInftyNonSep}, $(c2) \Rightarrow (c1)$. We prove 
$(c1) \Rightarrow (c2)$. Since $K \in {\cal NC}_k^*(G)$, by Claim \ref{NonSepToStar}, there exists $x \in V(G)$ such that $K= S(x, G)$. Since 
$K= S(x, G) \in {\cal C}_k^*(G)$, clearly $s(x, G) \le \rho^*_k(G) + 1$. However, if $s(x, G) = \rho^*_k(G) + 1$, then $K$ contains a cobase of $M_k(G)$, and so $E(G) \setminus K$ is a subset of a base of $M_k(G)$. Then 
$M_k(G \setminus x)$ has no $k$-circuit and therefore is not connected, contradicting $K \in {\cal NC}_k^*(G)$. Thus, $s(x, G) \le \rho^*_k(G)$.  

Since $K= S(x, G)$, we have: $M_k(G \setminus K) = M_k(G \setminus x)$. Since $K \in {\cal NC}_k^*(G)$, matroid $M_k(G \setminus x)$ is connected. 

Suppose that $k = 1$. Since $M_1(G \setminus x)$ is connected, by Theorem \ref{M1Con}, $G \setminus x \in {\cal CG}_{\bowtie}$. 

Now, suppose that $k \ge 2$. Since $M_k(G \setminus x)$ is connected, by Theorem \ref{MkCon}, $G \setminus x \in {\cal G}_{\bowtie}$. 
\ep

\vskip 1.5ex 

Let ${\cal S}_k(G)$ denote the set of vertex stars $S$ of $G$ such that
$|S| \le \rho_k^*(G)$.

\vskip 1.5ex 

It is known \cite{KplH} that ${\cal NC}^*(M(G)) = {\cal S}(G)$, where $M(G)$ is the cycle matroid of graph $G$. The next theorem is an analog of the above fact for $k$-circular matroids $M_k(G)$.

\begin{Theorem}
\label{nonsepeqstars}

{\em 
Let $G$ be  a $3$-connected graph, $k \geq 1$, and 
$|E(G)| - |V(G)| \ge k$. Then 
${\cal NC}^*(M_k(G)) = {\cal S}_k(G)$.
}

\end{Theorem} 

\bp (uses Theorems \ref{maxk}, \ref{MkCon}, \ref{M1Con}, and \ref{k>1NonSepStructure}) 

\vskip 0.3ex
Since $|E(G)| - |V(G)| \ge k$, by Theorem \ref{maxk}, $M_k(G)$ is a non-trivial matroid. Since $G$ is $3$-connected, clearly $G \in {\cal G}_{\bowtie}$. Then by Theorem \ref{MkCon} for $k \ge 2$ and Theorem \ref{M1Con} for $k = 1$, matroid $M_k(G)$ is connected. Therefore by Theorem \ref{k>1NonSepStructure}, ${\cal NC}^*(M_k(G)) \subseteq {\cal S}_k(G)$. 

\vskip 0.3ex
Since $G$ is $3$-connected, clearly $G \setminus x \in {\cal G}_{\bowtie}$ and 
$G \setminus x$ is a connected graph for every $x \in V(G)$. Therefore again by Theorem \ref{k>1NonSepStructure}, 
$ {\cal S}_k(G) \subseteq {\cal NC}^*(M_k(G))$. 

\vskip 0.3ex
Thus, ${\cal NC}^*(M_k(G)) = {\cal S}_k(G)$.  
\ep

\section{Uniquely representable matroids with 
${\cal S}(G) \subseteq {\cal NC}^*(G)$}
\label{UniquelyM1}

\indent 

Recall that a graph $G$ is uniquely defined by $M_k(G)$ if $M_k(G) = M_k(G')$ implies that graphs $G$ and $G'$ are strongly isomorphic.

\begin{claim}
\label{GuniqM}

{\em \cite{KplH}}
Let $G$ be a graph with  $|V(G)| \ge 4$. Suppose that $M(G)$ is a connected matroid. Then the following are equivalent:
\\[1ex]
$(c1)$ $G$ is uniquely defined by $M(G)$ and
\\[1ex]
$(c2)$ ${\cal S}(G) \subseteq {\cal NC}^*(G)$.

\end{claim}

In this part we will describe some results on matroids $M_k(G)$ analogous to the implication $(c2) \Rightarrow (c1)$ in Claim \ref{GuniqM} about matroid $M(G)$.

\vskip 1.5ex

We start with the following useful fact. 

\begin{claim}
\label{k>1S(x,G)inS(G')}

Let $G$ and $G'$ be graphs and  $k \ge 1$. Suppose that
\\[1ex]
$(a1)$ 
$M_k(G)$ is a connected matroid,   
\\[1ex]
$(a2)$ 
$x$ is a vertex of $G$ such that $S(x, G) \in {\cal NC}_k^*(G)$, and 
\\[1ex]
$(a3)$ $ M _k(G)= M _k(G')$.

\vskip 1ex
 Then 
$S(x, G) \in {\cal S}(G')$.

\end{claim}

\bp (uses Theorem \ref{k>1NonSepStructure} and Claim \ref{NonSepToStar}) 

Since $ M _k(G)= M _k(G')$, clearly 
${\cal NC}_k^*(G) ={\cal NC}_k^*(G')$. Since $S(x, G) \in {\cal NC}_k^*(G)$, we have: $S(x, G) \in {\cal NC}_k^*(G')$. Since $M_k(G')$ is a connected matroid, by Claim \ref{NonSepToStar}, there exists $x' \in V(G')$ such that $S(x, G) = S(x', G')$. Thus, $S(x, G) \in {\cal S}(G')$. 
\ep

\vskip 1.5ex

\begin{claim}
\label{eqnumver}
Let $G$ and $G'$ be graphs and  $k \ge 1$.   Suppose that
\\[1ex]
$(a1)$ $M_k(G)$ is a connected matroid and 
\\[1ex]
$(a2)$ $M_k(G) = M_k(G')$  $:$ $=$  $M_k$.

\vskip 1ex
Then  $|V(G)| = |V(G')|$.

\end{claim} 

\bp (uses Corollary \ref{grank})

Since $M_k(G) = M_k(G')$, we have: $\rho _k(G) = \rho _k(G')$. Since $M_k$ is a connected matroid, by Corollary \ref{grank}, $ \rho _k(G) = |V(G)| - 1 + k $ and  $\rho _k(G') = |V(G')| - 1 + k $. Therefore $|V(G)| = |V(G')|$.
\ep

\vskip 1.5ex 

Now, we are ready to prove our first results on graphs uniquely defined by their $k$-circular matroids analogous to implication $(c2) \Rightarrow (c1)$ in Claim \ref{GuniqM} about cycle matroid $M(G)$ of graph $G$. 

 \vskip 1ex

 We need the following known fact.
 
 \begin{claim}
\label{S(G)=S(G')}

{\em  \cite{KplH}}
Let $G$ and $G'$ be graphs and $E(G) = E(G')$. Then $G$ and $G'$ are strongly isomorphic if and only if ${\cal S}(G) = {\cal S}(G')$.

\end{claim}

\begin{Theorem} {\sc A condition for $M_k(G)$ to uniquely define graph $G$}
\label{AllNonSepUniq}

\vskip 0.3ex
Let $G$ and $G'$ be graphs and  $k \ge 1$. Suppose that
\\[1ex]
$(a1)$ 
$M_k(G)$ is a connected matroid,   
\\[1ex]
$(a2)$  
 ${\cal S}(G) \subseteq {\cal NC}_k^*(G)$,   and 
\\[1ex]
$(a3)$ $ M _k(G)= M _k(G')$.

\vskip 1ex
Then ${\cal S}(G) = {\cal S}(G')$, i.e. $G$ and $G'$ are strongly isomorphic.

\end{Theorem}

\bp (uses Claims \ref{S(G)=S(G')}, \ref{k>1S(x,G)inS(G')}, and \ref{eqnumver})   

Since $S(x, G) \in {\cal NC}_k^*(G)$ for every $x \in V(G)$,
by Claim \ref{k>1S(x,G)inS(G')},  ${\cal S}(G) \subseteq {\cal S}(G')$. Now, by Claim \ref{eqnumver}, $|V(G)| = |V(G')|$, and therefore $|{\cal S}(G)| = |{\cal S}(G')|$. It follows that ${\cal S}(G) = {\cal S}(G')$ and by Claim \ref{S(G)=S(G')}, $G$ and $G'$ are strongly isomorphic. 
 \ep

\vskip 1.5ex 

It is known \cite{W3} (see also \cite{KplH}) that multi-$3$-connected graphs $G$ are uniquely defined by $M(G)$, where $M(G)$ is the cycle matroid of graph $G$. The next theorem is an analog of the above fact for $k$-circular matroids $M_k(G)$.

\begin{Theorem} {\sc A condition for a graph $G$ to be uniquely defined by $M_k(G)$}
\label{k>1AllNonSepUniq}

\vskip 0.3ex
Let $G$ and $G'$ be graphs and  $k \ge 2$. Suppose that
 \\[1ex]
 $(a1)$ $\Delta(G) \ge k$ and $G \in {\cal G}_{\bowtie}$,
\\[0.5ex]
$(a2)$ $\Delta (G \setminus x) \ge k$ and $G \setminus x  \in {\cal G}_{\bowtie}$ for every vertex $x$ in $G$, 
  and 
  \\[0.5ex]
 $(a3)$
 $ M _k(G')= M _k(G)$.  
 
\vskip 0.6ex
 Then 
${\cal S}(G) = {\cal S}(G')$, i.e. $G$ and $G'$ are strongly isomorphic.

\end{Theorem} 

\bp (uses Theorems \ref{maxk}, \ref{MkCon},  \ref{k>1NonSepStructure}, and \ref{AllNonSepUniq}) 

By our assumption $(a1)$ and Theorems \ref{maxk} and \ref{MkCon}, we have: $M_k(G)$ is a connected matroid. Since $M_k(G)$ is a connected matroid, by Corollary \ref{grank}, we have: 
\vskip 0.7ex
$\rho ^*_k(G) = |E(G)| - |V(G)| + 1 - k$. 
\vskip 0.7ex
\noindent 
 Since 
$\Delta (G \setminus x) \ge k$ for every vertex $x$ in $G$, we have: 
\vskip 0.7ex
$  \Delta (G \setminus x) =  |E(G \setminus x)| - |V(G \setminus x)| = 
(|E(G)|  - s(x, G)) - (|V(G)| - 1) \ge k$ 
\vskip 0.7ex
\noindent 
for every vertex $x$ in $G$. Therefore  $ \rho^*_k(G) \ge s(x, G)$. By our assumption $(a2)$, $G \setminus x \in {\cal G}_{\bowtie}$ for every vertex $x$ in $G$. Since $s(x, G) \le \rho^*_k(G)$ for every vertex $x$ in $G$, we have by Theorem \ref{k>1NonSepStructure}:  
$S(x, G) \in {\cal NC}_k^*(G)$ for every $x \in V(G)$. Now, by Theorem \ref{k>1AllNonSepUniq}, ${\cal S}(G) = {\cal S}(G')$, i.e. $G$ and $G'$ are strongly isomorphic. 
\ep

\begin{Theorem} {\sc A condition for a graph $G$ to be uniquely defined by $M_1(G)$}
\label{k=1AllNonSepUniq}

\vskip 0.3ex
Let $G$ and $G'$ be graphs. Suppose that
 \\[1ex]
 $(a1)$ $\Delta(G) \ge 1$ and $G \in {\cal CG}_{\bowtie}$,
\\[0.5ex]
$(a2)$ $\Delta (G \setminus x) \ge 1$ and $G \setminus x  \in {\cal CG}_{\bowtie}$ for every vertex $x$ in $G$, 
  and 
  \\[0.5ex]
 $(a3)$
 $ M _1(G')= M _1(G)$.  
 
\vskip 0.6ex
 Then 
${\cal S}(G) = {\cal S}(G')$, i.e. $G$ and $G'$ are strongly isomorphic.

\end{Theorem} 

\bp (uses Theorems \ref{maxk}, \ref{M1Con},  \ref{k>1NonSepStructure}, and \ref{k>1AllNonSepUniq}) 

By our assumption $(a1)$ and Theorems \ref{maxk} and \ref{M1Con}, we have: $M_1(G)$ is a connected matroid. Since $M_1(G)$ is a connected matroid, by Corollary \ref{grank}, we have: $\rho ^*_k(G) = |E(G)| - |V(G)|$. Since 
$\Delta (G \setminus x) \ge 1$ for every vertex $x$ in $G$, we have: 
\vskip 0.7ex
$  \Delta (G \setminus x) =  |E(G \setminus x)| - |V(G \setminus x)| = 
(|E(G)|  - s(x, G)) - (|V(G)| - 1) \ge 1$ 
\vskip 0.7ex
\noindent 
for every vertex $x$ in $G$. Therefore  $ \rho^*_k(G) \ge s(x, G)$. By our assumption $(a2)$, $G \setminus x \in {\cal CG}_{\bowtie}$ for every vertex $x$ in $G$. Since $s(x, G) \le \rho^*_1(G)$ for every vertex $x$ in $G$, we have by Theorem \ref{k>1NonSepStructure}:  
$S(x, G) \in {\cal NC}_1^*(G)$ for every $x \in V(G)$. Now, by Theorem \ref{k>1AllNonSepUniq}, ${\cal S}(G) = {\cal S}(G')$, i.e. $G$ and $G'$ are strongly isomorphic. 
\ep

\vskip 1.5ex

\begin{Theorem}
\label{G3-conAllSmallUniq}
 
 Let $G$ and $G'$ be graphs and  $k \ge 1$. 

Suppose that
 \\
 $(a1)$ $G$ is a $3$-connected graph,
\\[0.5ex]
$(a2)$ ${\cal C}_k(G) \ne \emptyset $, 
 \\[0.5ex]
 $(a3)$ $s(x, G) \le \rho^*_k(G)$ for every vertex $x$ of $G$, and 
 \\[0.5ex]
 $(a4)$ $M _k(G')= M _k(G)$. 
 
\vskip 0.6ex
  Then 
${\cal S}(G) = {\cal S}(G')$, i.e. $G$ and $G'$ are strongly isomorphic.

\end{Theorem}

\bp (uses Theorems \ref{maxk}, \ref{MkCon}, \ref{M1Con}, \ref{k>1AllNonSepUniq}, and \ref{k=1AllNonSepUniq}) 

Since $G$ is a $3$-connected graph, $G \in {\cal CG}_{\bowtie}$. Since 
${\cal C}_k(G) \ne \emptyset $, by Theorem \ref{maxk}, matroid $M_k(G)$ is non-trivial. Therefore, by Theorems \ref{MkCon} for $k \ge 2$ and \ref{M1Con} for $k = 1$, $M_k(G)$ is a connected matroid. Clearly, all assumptions of Theorems \ref{k>1AllNonSepUniq} and  \ref{k=1AllNonSepUniq}  hold. Thus, ${\cal S}(G) = {\cal S}(G')$, i.e. $G$ and $G'$ are strongly isomorphic.   
\ep

\begin{remark}
\label{AllNonsepUniqG}

{\em
Theorems \ref{k>1AllNonSepUniq} - \ref{G3-conAllSmallUniq} describe various conditions that guarantee a unique graph representation of some $k$-circular matroids.  
}

\end{remark}

\section{Uniquely representable matroids with 
${\cal S}(G) \not \subseteq {\cal NC}^*(G)$}
\label{UniquelyM2}

\indent

Now, we will describe some results on matroids $M_k(G)$ illustrating a new phenomenon when implication 
$(c1) \Rightarrow (c2)$ in Claim \ref{GuniqM} does  not hold if matroid $M(G)$ is replaced by $M_k(G)$.

\vskip 1.5ex

\begin{Theorem}
\label{G,Fnoloops}

Let $G$ and $G'$ be graphs and  $k \ge 1$. 

\vskip 0.3ex
Suppose that
 \\
 $(a1)$ 
 both graphs $G$ and $G'$ have no loops, 
 \\[0.5ex]
 $(a2)$
  $M_k(G)$ is a connected matroid, 
  \\[0.5ex]
 $(a3)$ $G$ has exactly one vertex $v$ such that 
 $S(v, G) \notin {\cal NC}^*_k(G)$, and 
  \\[0.5ex]
 $(a4)$
 $ M _{k}(G')= M _k(G)$.
\vskip 0.6ex
 Then 
${\cal S}(G) = {\cal S}(G')$, i.e. $G$ and $G'$ are strongly isomorphic.

\end{Theorem}

\bp (uses Claims  \ref{S(G)=S(G')}, \ref{k>1S(x,G)inS(G')}, and \ref{eqnumver})

If $x$ is a vertex of $G$ distinct from $v$, then by Claim \ref{k>1S(x,G)inS(G')}, $S(x, G) \in {\cal S}(G')$. Then 
${\cal S}(G) \setminus \{S(v, G)\}  \subseteq {\cal S}(G')$. 
 By Claim \ref{eqnumver}, $|V(G)| = |V(G')|$, and therefore 
 $|{\cal S}(G)| = |{\cal S}(G')|$. Therefore only one vertex, say $v'$, in $G'$ is such that  $S(v', G') \not \in {\cal S}(G) \setminus \{S(v, G)\}$. Then 
 ${\cal S}(G') \setminus \{S(v', G')\} = {\cal S}(G) \setminus \{S(v, G)\}$.
 Thus, it remains to prove that $S(v, G) = S(v', G')$.   
        
 \vskip 1ex
 
Let $e \in S(v', G')$. Then clearly, $e$ belongs to at most one element of 
${\cal S}(G') \setminus\{S(v', G')\}$.
 Since 
${\cal S}(G') \setminus \{S(v', G')\} = {\cal S}(G) \setminus \{S(v, G)\}$, we have: $e$ belongs to at most one element of ${\cal S}(G) \setminus \{S(b, G )\}$. Since $G$ has no loops, $e$ is incident to $b$.  Hence $S(v', G') \subseteq S(v, G) $. 
 
  \vskip 1ex
    
Since both $G$ and $G'$ have no loops, by the arguments similar to those above, 
\\
$S(v, G) \subseteq S(v', G') $. Thus, $S(v, G) = S(v', G') $. Hence ${\cal S}(G) = {\cal S}(G')$, and by Claim \ref{S(G)=S(G')}, $G$ and $G'$ are strongly isomorphic.
\ep

\begin{Theorem}
\label{AlmostAllNonsep} {\sc Another condition for $M_k(G)$ to uniquely define graph $G$}

\vskip 0.7ex
Let $G$ and $G'$ be graphs and  $k \ge 1$. Suppose that
 \\[0.5ex]
 $(a1)$
  $M_k(G)$ is a connected matroid, 
  \\[0.5ex]
 $(a2)$ $S(x, G) \in {\cal NC}^*_k(G)$ for every vertex $x$ of $G$ except for one vertex $v$,
 \\[0.5ex]
 $(a3)$ if $p$ is a loop in $G$, then $p$ is incident to $v$, 
 \\[0.5ex]
 $(a4)$ $v$ is not a $k$-big vertex of $G$ and $G \setminus v$ has no tree component, and 
  \\[0.5ex]
 $(a5)$
 $ M _{k}(G')= M _k(G)$.
\vskip 0.6ex
 Then 
${\cal S}(G) = {\cal S}(G')$, i.e. $G$ and $G'$ are strongly isomorphic.

\end{Theorem} 

\bp (uses Theorems \ref{basstructure}, \ref{SmallTightCo} and \ref{G,Fnoloops} and Claims \ref{LeafCocircuits}, \ref{k>1S(x,G)inS(G')}, and \ref{eqnumver})

Our first step in proving our theorem is to prove that there exists a bijection 
\\
$\alpha:V(G) \to V(G')$ such that $S(x, G) = S(\alpha (x), G')$ for every $x \in V(G) \setminus v$.

If $x$ is a vertex of $G$ distinct from $v$, then by Claim \ref{k>1S(x,G)inS(G')}, $S(x, G) \in {\cal S}(G')$, and so 
${\cal S}(G) \setminus \{S(v, G)\}  \subseteq {\cal S}(G')$. Therefore for $x \in V(G) \setminus v$ there exists a unique vertex $x'$ such that 
$S(x', G') = S(x, G)$. In this case we put $\alpha (x) = x'$. Since $M_k(G)$ is a connected matroid and $ M _{k}(G')= M _k(G)$, by Claim \ref{eqnumver}, 
$|V(G)| = |V(G')|$, and therefore $|{\cal S}(G)| = |{\cal S}(G')|$. Hence only one vertex, say $v'$, in $G'$ is such that  
$S(v', G') \not \in {\cal S}(G) \setminus \{S(v, G)\}$. Clearly,  
$v' \ne \alpha (x)= x'$ for every $x \in V(G) \setminus v$. We put 
 $v' = \alpha (v)$. Clearly, $\alpha$ is a bijection from $V(G)$ to $V(G')$ and
 $S(x', G') = S(x, G)$ for every $x \in V(G) \setminus v$.  
 
  \vskip 1ex
  
Now, we will prove that $\alpha$ is a strong isomorphism from $G$ to $G'$. Obviously, it is sufficient to prove that 
$S(v, G) = S(\alpha (v), G') = S(v', G')$. 
\\[1ex]
${\bf (p1)}$ First, suppose that $G$ has no loops. By Theorem \ref{G,Fnoloops}, it is sufficient to prove that $G'$ has no loops. Suppose not. Then $G'$ has a loop, say $e$. Since $G$ has no loops, $e$ is not a loop in 
$G$, and so $e$ is incident to two distinct vertices in $G$, say $x_1$ and $x_2$. We claim that $v \in \{x_1, x_2\}$. Indeed, suppose not. Then  
$v \notin \{x_1, x_2\}$. Therefore $e \in S(x_1, G) \cap S(x_2, G)$, and so 
$e \in S(x'_1, G') \cap S(x'_2, G')$. Since 
 $x'_1 = \alpha (x_1) \ne \alpha (x_2) = x'_2$, $e$ is not a loop in $G'$, a contradiction. Thus, $e$ in $G$ is incident to $v$ and to another vertex of $G$, say $z$. 
 
Since $e \in S(z, G) = S(z', G')$ and $e$ is a loop in $G'$, clearly, edge $e$ is a loop at $z' = \alpha (z)$ in $G'$. 
Now, by assumption $(a4)$ of our theorem, $|S(v, G)| \le \rho_k^*(G)+1$ and 
$G \setminus v$ has no tree component. By Theorem  \ref{SmallTightCo}, 
$S(v, G) \in {\cal C}^*_k(G)$. Therefore, by Claim \ref{LeafCocircuits}, there exists $B \in {\cal B}_k(G)$ such that $e$ is a dangling edge at $v$ in 
$G \langle B \rangle$. Since $ M _{k}(G')= M _k(G)$, clearly,
 $B \in {\cal B}_k(G)$.  Since $G$ has no loops and $S(x', G') = S(x, G)$ for every $x \in V(G) \setminus v$, we have: 
$G \langle B \rangle \setminus v = G' \langle B \rangle \setminus v'$. Therefore $v'$ is an isolated vertex in $G' \langle B \rangle$, and so $B$ does not span $V(G')$, contradicting Theorem \ref{basstructure}. 
\\[1ex]
${\bf (p2)}$ Now, suppose that $G$ has a loop. Since $S(x', G') = S(x, G)$ for every 
$x \in V(G) \setminus v$, an edge $p$ is a loop at $v$ in $G$ if and only if $p$ is a loop at $v'$ in $G'$. 

First, suppose that $e \in S(v, G)$ and $e$ is not a loop in $G$.  By our assumption $(a4)$ and Theorem \ref{SmallTightCo}, 
$S(v, G) \in {\cal C}^*_k(G)$. Therefore, by Theorem \ref{LeafCocircuits}, there exists $B \in {\cal B}_k(G)$ such that $e$ is a dangling edge at $v$ in $G \langle B \rangle$. Clearly, from our assumption $(a3)$ and  equality 
$S(x', G') = S(x, G)$ for every $x \in V(G) \setminus v$ we have: 
$G \langle B \rangle \setminus v = G' \langle B \rangle \setminus v'$. If $e$ is not incident to $v'$ in $G'$, then $B$ does not span $V(G')$, contradicting Theorem \ref{basstructure}. Thus, $e \in S(v', G')$.  

Finally, suppose that $e \in S(v', G')$ and $e$ is not a loop in $G'$. Since $S(x', G') = S(x, G)$ for every 
$x \in V(G) \setminus v$, clearly, $e$ belongs to at most one element of 
${\cal S}(G) \setminus\{S(v, G)\}$. By our assumption $(a3)$, edge $e$ is incident to $v$ in $G$. Therefore $e \in S(v, G)$. 
 \ep

\vskip 2ex

Using Theorem \ref{AlmostAllNonsep} and the arguments similar to those in the proof of Theorem \ref{k>1AllNonSepUniq}, one can prove  the following theorem.

\begin{Theorem}
\label{k>1AlmostAllNonsepG} {\sc Another condition for $G$ to be uniquely defined by $M_k(G)$}

\vskip 0.7ex
Let $G$ and $G'$ be graphs and  $k \ge 2$. Suppose that
 \\[0.5ex]
 $(a1)$
  $\Delta(G) \ge k$ and $G \in {\cal G}_{\bowtie}$, 
  \\[0.5ex]
 $(a2)$ $\Delta (G \setminus x) \ge k$ and $G \setminus x  \in {\cal G}_{\bowtie}$ for every vertex $x$ in $G$ except for one vertex $v$,
 \\[0.5ex]
 $(a3)$ if $p$ is a loop in $G$, then $p$ is incident to $v$, 
 \\[0.5ex]
 $(a4)$ $v$ is not a $k$-big vertex of $G$ and $G \setminus v$ has no tree component, and 
  \\[0.5ex]
 $(a5)$
 $ M _{k}(G')= M _k(G)$.
\vskip 0.6ex
 Then 
${\cal S}(G) = {\cal S}(G')$, i.e. $G$ and $G'$ are strongly isomorphic.

\end{Theorem} 

\vskip 0.1ex

Using Theorem \ref{AlmostAllNonsep} and the arguments similar to those in the proof of Theorem \ref{k=1AllNonSepUniq}, one can prove  the following theorem.

\begin{Theorem}
\label{k=1AlmostAllNonsepG} {\sc Another condition for $G$ to be uniquely defined by $M_1(G)$}

\vskip 0.7ex
Let $G$ and $G'$ be graphs. Suppose that
 \\[0.5ex]
 $(a1)$
  $\Delta(G) \ge 1$ and $G \in {\cal CG}_{\bowtie}$, 
  \\[0.5ex]
 $(a2)$ $\Delta (G \setminus x) \ge 1$ and $G \setminus x  \in {\cal CG}_{\bowtie}$ for every vertex $x$ in $G$ except for one vertex $v$,
 \\[0.5ex]
 $(a3)$ if $p$ is a loop in $G$, then $p$ is incident to $v$, 
 \\[0.5ex]
 $(a4)$ $v$ is not a $1$-big vertex of $G$ and $G \setminus v$ has no tree component, and 
  \\[0.5ex]
 $(a5)$
 $ M _{1}(G')= M _1(G)$.
\vskip 0.6ex
 Then 
${\cal S}(G) = {\cal S}(G')$, i.e. $G$ and $G'$ are strongly isomorphic.

\end{Theorem}

Using Theorems  \ref{maxk}, \ref{MkCon}, \ref{M1Con}, \ref{k>1AlmostAllNonsepG}, and \ref{k=1AlmostAllNonsepG} and the arguments similar to those in the proof of Theorem \ref{G3-conAllSmallUniq}, one can prove  the following theorem.

 \begin{Theorem}
\label{G3-conAlmostAllNonsep}

Let $G$ and $G'$ be graphs, $t \in V(G)$, and $k \ge 1$.
Suppose that
\vskip 0.3ex
\noindent
 $(a1)$ 
 graph $G$ is $3$-connected, 
 \\[0.5ex]
 $(a2)$
  ${\cal C}_k(G) \ne \emptyset $, 
  \\[0.5ex]
 $(a3)$
 $s(t, G) = \rho^*_k(G) + 1$
  and  $s(x, G) \le \rho^*_k(G)$ for every vertex $x$ of $G$ distinct from $t$, and 
  \\[0.5ex]
 $(a4)$
 $ M _{k}(G')= M _k(G)$.
\vskip 0.6ex
 Then 
${\cal S}(G) = {\cal S}(G')$, i.e. $G$ and $G'$ are strongly isomorphic.

\end{Theorem}

\begin{claim}
\label{wheel}

{\em (2.3.2 in \cite{KK1})}

\vskip 0.3ex
Let  $G$ be a graph with $v(G) = v$ and $e(G) = e$.  
Suppose that $G$ is 3-connected  and has a vertex $x$ such that 
$d(x,G) > e - v$.
Then  $G$ is the wheel with center $x$, and so $d(x,G) = v - 1$ and 
$d(z,G) = 3$ for every $z \in V(G) \setminus x$.

\end{claim}

From Claim \ref{wheel} and Theorems \ref{G3-conAllSmallUniq} and \ref{G3-conAlmostAllNonsep} we have, in particular:

\begin{Theorem}
\label{3-ConBiciUni}

Every 3-connected graph $G$ not isomorphic to a complete graph on 4 vertices is uniquely defined by its bicircular matroid $M_1(G)$.

\end{Theorem}

Another proof of Theorem \ref{3-ConBiciUni} can be obtained using some results in  \cite{CGW}.

\begin{figure}[hptb]
\begin{center}
\scalebox{0.35}[.35]{\includegraphics{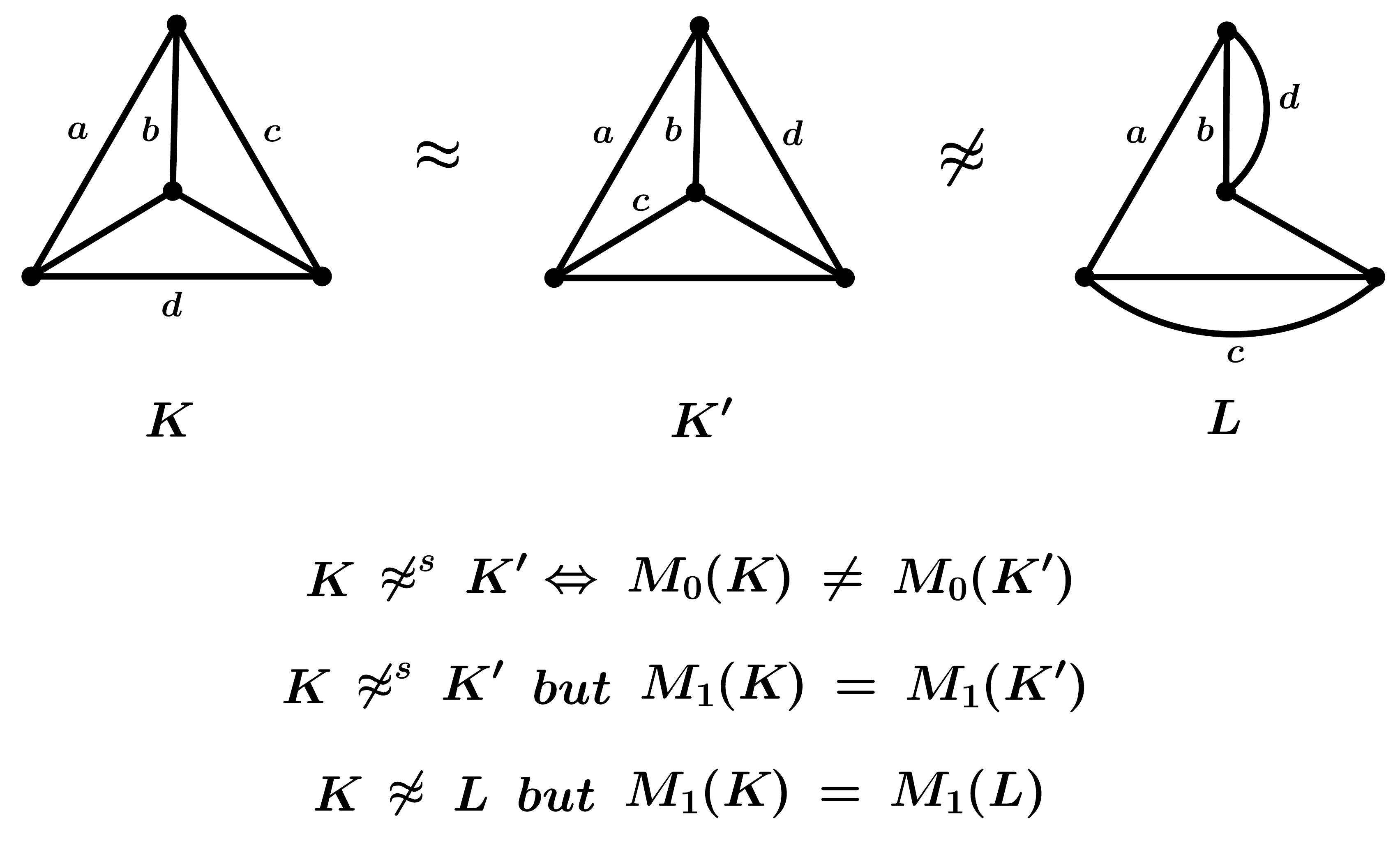}}
\end{center}
\caption
{Graph $K$ is not uniquely defined by matroid $M_1(K)$.}
\label{K4fork=0,1Fig}
\end{figure}

\begin{Theorem} {\sc When is a 3-connected graph $G$ uniquely defined by $M_k(G)$ ?}
\label{mainth}

\vskip 0.3ex
Let $G$ be a graph, $k \ge 1$, and $M_k(G)$ the $k$-circular matroid of $G$.
Suppose that $G$ is 3-connected. 
 Then the following are equivalent:
\\[0.7ex]
$(c1)$  $G$ is uniquely defined by  $M_k(G)$ and
\\[0.7ex]
$(c2)$ $|E(G) | - |V(G)| \ge k$ and either every vertex of $G$ is $k$-small or every vertex of $G$ is $k$-small except for one which is $k$-tight. 
  
\end{Theorem}

The implication 
$(c2) \Rightarrow (c1)$ of Theorem \ref{mainth} follows immediately from Theorems 
\ref{G3-conAllSmallUniq} and \ref{G3-conAlmostAllNonsep}. In our next paper we will describe some graph operations that will provide non-isomorphic graphs with the same $k$-circular matroid. Implication 
$(c1) \Rightarrow (c2)$ in Theorem \ref{mainth} will follow from those results.

\vskip 1.5ex
Theorem \ref{mainth} is a natural extension of the classical Whitney's matroid-isomorphism theorem on the cycle matroid of a $3$-connected graph \cite{W3} (see also\cite{KplH}).


\addcontentsline{toc}{chapter}{Bibliography}


\begin{thebibliography}{99}


\bibitem{BM} J. A. Bondy and U. S. R. Murty, \emph{Graph Theory},
3rd Corrected Printing, GTM 244, Springer-Verlag,
New York, 2008.

\bibitem{CGW} C. R. Coulard, J. G. del Greco, D. K. Wagner, 
Representations of bicircular matroids,
{\em  Discrete Applied Mathematics} {\bf 32}
(1991) 223-240.

\bibitem{D} R. Diestel, {\em Graph Theory}, Springer-Verlag,
New York, 2000.

\bibitem{HalJ} R. Halin and H. A.  Jung,  Note on isomorphisms of graphs,
{\em J. London Math.} Sec. {\bf 42} (1967) 254-256.

\bibitem{HJK}  R. L. Hemminger, H. A.  Jung, A. K. Kelmans,
 On 3--skein isomorphisms of a graph, 
{\em Combinatorica} {\bf 2} (4), (1982) 373-376.

\bibitem{KK1} J. F. De Jes\'us and A. K. Kelmans, $K$-circular matroids of graphs, arXiv:1508.05364 (2015). 

\bibitem{Ksisd} 
A. K. Kelmans, On semi-isomorphisms
and semi-dualities of graphs,  
%
{\em Graphs and
Combinatorics} {\bf 10} (1994) 337--352. 

\bibitem{K3skHng} A. K. Kelmans, On 3-skeins in a 3-connected graph. 
%
{\em  Studia Scientarum
Mathematicarum Hungarica} {\bf 22} (1987) 265-273.

\bibitem{KplH} 
A. K. Kelmans, The concept of a vertex
in a matroid, the non-separating cycles  and a new
criterion for graph planarity. 
%\\
In {\em  Algebraic
Methods in Graph Theory},
 Vol. 1, Colloq. Math. Soc. J\'anos Bolyai, (Szeged, Hungary, 1978)   
North--Holland {\bf 25} (1981) 345-388.

 
\bibitem{Ox} J. G. Oxley, {\em Matroid Theory}, Oxford University Press, 2006. 

\bibitem{Wag} D. K. Wagner, Connectivity  in bicircular matroids,
{\em J. Combinatorial Theory} {\bf B 39} (1985) 308-324.

\bibitem{Welsh} D. J. A.  Welsh, {\em Matroid Theory},
Academic Press, London, 1976.


\bibitem{W1}   
H. Whitney, On the abstract properties of linear dependence,
{\em Amer. J. Math.} {\bf 57} (1935) 509-533.

\bibitem{W2}  
H. Whitney, 2-isomorphic graphs,
{\em Amer. J. Math.} {\bf 55} (1933) 245-254.

\bibitem{W3}  
H. Whitney, Congruent graphs and the connectivity of graphs,
{\em Amer. J. Math.} {\bf 54} (1932) 159-168.

\bibitem{W4}  
H. Whitney, Non-separable and planar graphs,
{\em Amer. J. Math.} {\bf 34} (1932) 339-362.

\end{thebibliography}
\end{document}